\documentclass[12pt]{amsart}

\usepackage{enumitem} 
\usepackage{xspace} 
\usepackage{listings} 

\newcommand \GBasis {{Gr\"obner basis}\xspace}
\newcommand \GBases {{Gr\"obner bases}\xspace}
\newcommand\cocoa{{CoCoA}\xspace}
\newcommand\cocoav{{CoCoA-5}\xspace}
\newcommand\cocoalib{{CoCoALib}\xspace}
\newcommand\cocoacode[1]{\protect{``\texttt{#1}''}}
\newcommand \ie {\textit{i.e.}}
\newcommand \eg {\textit{e.g.}}
\newcommand \ideal[1] {\langle #1 \rangle}
\newcommand \QQ {{\mathbb Q}}
\newcommand \ZZ {{\mathbb Z}}
\newcommand \LPP {\mathop{\rm LT}\nolimits}
\newcommand{\gin}{{\mathrm{gin}}}
\newcommand{\minpoly}[1]{{\mu_{#1}}}

\usepackage{amssymb}
\usepackage{amsmath}

\title[Gr\"obner Bases in CoCoA and CoCoALib]
{Gr\"obner Bases for Everyone with CoCoA-5 and CoCoALib}

\author[J. Abbott, A.M. Bigatti]{John Abbott, Anna Maria Bigatti}

\address{John Abbott\\Institut f\"ur Mathematik\\Universit\"at Kassel
}

\email{abbott@dima.unige.it}

\address{Anna Maria Bigatti\\Dipartimento di Matematica\\Universit\`a
  degli Studi di Genova
}

\email{bigatti@dima.unige.it}

\subjclass[2010]{
13-04, 
13P10, 
68W30, 
13P15  
14Q10, 
05E40, 
68W25 
}

\keywords{Gr\"obner bases, elimination, term orderings, modular
  methods, software library}

\thanks{(J.~Abbott) J.~Abbott is an INdAM-COFUND Marie Curie Fellow}
\thanks{(Anna M. Bigatti) This work was partly
  supported by the ``National Group for Algebraic and Geometric
  Structures, and their Applications'' (GNSAGA -- INdAM)}

\begin{document}

\begin{abstract}
We present a survey on the developments related to \GBases, and show explicit
examples in \cocoa.

The CoCoA project dates back to 1987: its aim was to create a
``mathematician''-friendly computational laboratory for
studying Commutative Algebra, most especially Gr\"obner bases.
Always maintaining this ``friendly'' tradition, the project has grown
and evolved, and the software has been completely rewritten.

\cocoa offers
Gr\"obner bases for all levels of interest: from the basic, explicit call in
the interactive system CoCoA-5~\cite{CoCoA-5}, to problem-specific
optimized implementations, to the computer--computer communication
with the open source C++ software library, CoCoALib~\cite{CoCoALib}, or
the prototype OpenMath-based server.

The openness and clean design of CoCoALib and CoCoA-5 are intended to
offer different levels of usage, and to encourage external contributions.
\end{abstract}

\maketitle

\section{Introduction}

The CoCoA project traces its origins back to 1987 under the lead of L.~Robbiano:
the aim was to create a software laboratory for studying Commutative
Algebra and especially Gr\"obner bases, and which is welcoming even to
mathematicians who are wary of new-fangled computers.

Since then the realm of applicability of \GBases has
continually expanded, so researchers interested in using them
now come from a broad palette of subject areas ranging from
the theoretical to quite practical topics.
So there are still
``pure'' mathematicians as at the outset, but now also ``programming''
mathematicians, and statisticians, computer scientists, and so on.
Another factor crucial in making \GBases relevant to practical problems is
the interim progress in computer hardware and software techniques.

Since its beginning, the \cocoa project has evolved considerably, and the software has been
rewritten: it now comes in the form of the very flexible software combination
\cocoav/\cocoalib, while maintaining its tradition of being
\textit{user-friendly} so it offers \GBases for all levels of
interest and programming ability~---~\textit{a \GBasis for everyone}.
This means that the ``\cocoa experience'' covers a
wide range: from the basic, explicit call
in the interactive system CoCoA-5~\cite{CoCoA-5}
(see Section~\ref{sec:GBwithEase}),
to functions which use \GBases implicitly
(see Sections~\ref{sec:gfan} and~\ref{sec:gin}),
to problem-specific optimized implementations
(see Sections~\ref{sec:elim}, \ref{sec:BM} and~\ref{sec:rat}),
to the computer--computer
communication with the open source C++ software library,
CoCoALib~\cite{CoCoALib}, or with the prototype OpenMath-based server
(see Section~\ref{sec:cocoalib}).

The importance that \GBases have acquired derives from the fact they
enable or facilitate so many other computational mathematical results.
A natural consequence is that a \GBasis is \emph{almost never the final
answer} that is sought, but just a stepping stone on the way to the goal,
\eg~a Hilbert series or a primary decomposition.
In this paper we concentrate on those computations in \cocoa which are
directly related to Gr\"obner bases, illustrating the wide range of
applications which have evolved over the last 50 years, and
providing (explicitly or implicitly) Gr\"obner bases for everyone.

\subsection{What is new in \cocoav?  And what is not?}

CoCoA-4 was widely appreciated for its ease of use, and the naturalness of its
interactive
language.  However, it did have limitations, and several ``grey
areas''.
We designed the new CoCoA-5 language to strike a
balance between backward-compatibility (to avoid alienating existing CoCoA-4 users)
and greater expressibility with a richer and more solid mathematical
basis (eliminating those ``grey areas'').

So, what's not new?  Superficially the new CoCoA-5 language and system
closely resemble CoCoA-4 because we kept it largely backward
compatible.  At the same time CoCoA-5 improves upon the
underlying mathematical structure and robustness
of the old system.  We are very aware that a number of CoCoA users
are mathematicians with only limited programming experience, for whom
learning CoCoA was a ``big investment'', and who are reluctant to make
another such investment~---~that is why we wanted to make the passage
to \cocoav as painless as possible.

So, if almost nothing has changed, \emph{what's new}?
The clearly defined semantics of the \cocoav
new language make it both more robust and more flexible;
it provides greater expressibility and a more solid mathematical basis.
In particular, it offers full flexibility for the field of
coefficients: \eg~$\ZZ/\ideal{p}$ with large~$p$,
fraction fields and algebraic extensions
(see Section~\ref{sec:morecoefficients}),
and even heuristically verified floating point arithmetics
with rational reconstruction (see Section~\ref{sec:twinfloats}).

However, under the surface, the change is radical since its
mathematical
core, \cocoalib, has been rewritten from scratch, to be
faster, cleaner and more powerful than the old system, and
also to be used as a C++ library.

\subsection{How Do CoCoALib and CoCoA-5 differ?}
\label{sec:cocoalib-vs-cocoa5}

The glib answer is: \textit{As little as possible!}

\medskip


One important idea behind the designs of \cocoalib and \cocoav is
  that of making it easy to take a prototype implementation in CoCoA-5 and
  translate it into C++ using CoCoALib.  
We intend to facilitate the ``translation'' step as much as possible.

\cocoalib~\cite{CoCoALib,Abb06}, the C++ library,
contains practically all the mathematical knowledge
and ability whereas \cocoav~\cite{CoCoA-5}
offers convenient, interactive access to \cocoalib's capabilities.
Most functions are accessible from both, and
have identical names and behaviour
(see Section~\ref{sec:cocoalib}).

More precisely,
\cocoav is an interactive, interpreted environment which makes it
better suited to ``rapid prototyping'' than the relatively rigid, statically typed regime
of C++.
To keep it simple to learn, \cocoav has only a few data types: for instance,
a power-product in \cocoav is represented as a monic polynomial with a single term,
\ie~a ring element (of a polynomial ring).  In contrast,
in \cocoalib there is a dedicated class, \texttt{PPMonoidElem}, which
directly represents each power-product, and allows efficient
operations on the values (\eg~without the overhead of the superfluous
coefficients present in the simplistic approach of CoCoA-5).

Programming with \cocoalib does tend to be more onerous than with \cocoav, largely
because of C++'s demanding, rigid rules.  However, the reward is
greater flexibility and typically faster
computation (sometimes much faster).  Also, of course, those who want to use
\cocoalib's abilities in their own C++ program necessarily have to use \cocoalib.
This is why our goal is that everything which can be computed in \cocoav
should be just as readily computable with CoCoALib.
Currently a few \cocoav functions are still implemented in
\cocoav packages, but these are being steadily translated into C++.
(We don't say which ones, because the list is constantly shrinking)

\section{Gr\"obner Bases with Ease}
\label{sec:GBwithEase}

The simplest cases of Gr\"obner bases are
for ideals in $\QQ[x_1,\ldots,x_n]$ or $\ZZ/\ideal{p}[x_1,\ldots,x_n]$ with~$p$ a prime.
These are also the easiest cases to give to \cocoa.
Here is an example:  
\begin{lstlisting}[alsolanguage=C++]
/**/ use QQ[x,y,z];   // or ZZ/(2)[x,y,z];
/**/ I := ideal(x^3 +3,  y-x^2,  z-x-y);
/**/ GBasis(I);
[x +y -z,  y^2 -3*y +3*z,  y*z -3*y +3*z +3,  z^2 -4*y +3*z +6]
\end{lstlisting}

An essential ingredient in the definition of a \GBasis is the
term-ordering: \ie~a total ordering on the power-products which respects
multiplication, and where 1 is the smallest power-product.  

In this example the term-ordering was not explicitly indicated, so CoCoA assumes
\texttt{StdDegRevLex} 
(with the common convention that the 
indeterminates generating the polynomial ring were given in decreasing
order).
In many computer algebra systems this is the
default ordering because it 
\textit{generally} gives the best performance and most compact answer.
Another well-known family of orderings is
\texttt{lex} (short for ``lexicographic'').  
A \texttt{lex} \GBasis 
of a zero dimensional ideal in normal position
has a particular \textit{shape}
which is theoretically useful for solving polynomial systems
(see, for example, the Kreuzer--Robbiano
book~\cite{KR00}, Sec.~3.7).  However its practical usefulness is limited by
the fact that \texttt{lex} bases tend to be particularly big and ugly,
and are frequently rather costly to compute.

There are other gradings and orderings
which are useful for studying specific problems:
for instance, an important family
are the \textit{elimination orderings} which are used implicitly
in Section~\ref{sec:elim}.
\cocoa also offers a fully general,
matrix-based implementation of term-orderings
(see Sections~\ref{sec:gfan} and~\ref{sec:elim}).

\medskip

In \cocoa the term-ordering 
is specified at the same time as the polynomial ring; \GBases of all
ideals in that polynomial ring will automatically 
be computed with respect to that ordering.
Thus, in \cocoa the term-ordering is an intrinsic property of each polynomial
ring.  This means that $\QQ[x,y,z]$ with \texttt{lex} is regarded as a different
ring from $\QQ[x,y,z]$ with \texttt{StdDegRevLex}.
Here is an example of computing a \texttt{lex} \GBasis.

\goodbreak
\begin{lstlisting}[alsolanguage=C++]
/**/ use QQ[x,y,z], lex; // specify ordering together with ring
/**/ I := ideal(x^3 +3,  y-x^2,  z-x-y);
/**/ ReducedGBasis(I);   // basis is wrt lex ordering
[x +(1/4)*z^2 +(-1/4)*z +3/2,  y +(-1/4)*z^2 +(-3/4)*z -3/2,
   z^3 +9*z -6]
\end{lstlisting}
In the last example above we used the command \texttt{ReducedGBasis} which
computes a reduced \GBasis: namely a ``cleaned up'' basis
with only non-redundant, monic, fully reduced
elements~---~it is unique (up to the order of its elements).

\subsection{Verbosity and interruption}\label{sec:verbosity}

Sometimes it is handy to know what is happening inside a running function.
For example, a Gr\"obner basis computation may be taking a long time,
and we would like to know whether it is likely far from finishing, and if
so, interrupt it.

A new feature in CoCoA-5.2.0 is the ability to set the verbosity level;
there is another function which tells you the current level.

\begin{lstlisting}[alsolanguage=C++]
/**/ SetVerbosityLevel(100);
/**/ VerbosityLevel();
100
\end{lstlisting}

This is a global setting, and higher verbosity levels trigger the
printing of increasing amounts of internal ``progress information'' in several
functions (both in CoCoA-5 and in CoCoALib).

For instance, the lowest level giving information on the progress of \GBases
is~100; everytime a new polynomial is found,
a line like this is printed:
\begin{lstlisting}[alsolanguage=C++]
myDoGBasis[1]: New poly in GB: len(GB) = 10 len(pairs) = 6
\end{lstlisting}

In some hard \GBasis computations, by setting the verbosity level,
one may see that the number of pairs yet to
be processed is unfeasibly high.  
The user may then choose to interrupt the computation
by typing \texttt{Ctrl-C}:
the computation will be interrupted as soon as the reduction of the
current S-polynomial terminates.


This interruption cancels the incomplete \GBasis computation,
and returns the computer to the state it was in just before the
\GBasis computation was begun
(thanks to the clean, exception-safe design of CoCoALib).

Besides \texttt{GBasis}, verbose information can be produced by
numerous functions: see, for instance, Sections~\ref{sec:gfan}
and~\ref{sec:gin}.  Indeed, the number of functions (both in \cocoalib
and \cocoav) which respond to the verbosity setting is steadily
increasing~---~details are in the documentation (type ``\texttt{?verbose}'').
Similarly the number
of interruptible \cocoalib functions is gradually increasing; in any case, all
interpreted \cocoav functions can be interrupted.

\subsection{More rings of coefficients}
\label{sec:morecoefficients}

The easy examples above show the definition of a polynomial rings
with rational coefficients, but the choice of coefficients in \cocoa is quite wide.
For example, coefficients in a finite field $\ZZ/\ideal{p}$:
\begin{lstlisting}[alsolanguage=C++]
/**/ use ZZ/(10^29 + 319)[x];
/**/ ReducedGBasis(ideal(3*x-1));
[x -33333333333333333333333333440]
\end{lstlisting}

\noindent
Or coefficients in algebraic extension fields:
\begin{lstlisting}[alsolanguage=C++]
/**/ use R ::= QQ[i];
/**/ K := R/ideal(i^2 +1);
/**/ use K[x,y,z];
/**/ I := ideal(i*x^3 -z, x^2*y^3 -i*y*z^2);
/**/ ReducedGBasis(I);
[x^3 +(i)*z,  x^2*y^3 +(-i)*y*z^2,  y^3*z +x*y*z^2]

/**/ use R ::= QQ[sqrt2, sqrt3];
/**/ K := R/ideal(sqrt2^2 -2, sqrt3^2 -3);
/**/ IsField(K);
true
/**/ use K[x,y,z];
/**/ I := ideal(sqrt3*x^2 -y, x*y -sqrt2*z);
/**/ ReducedGBasis(I);
[x*y +(-sqrt2)*z,  x^2 +((-1/3)*sqrt3)*y,  y^2 +(-sqrt2*sqrt3)*x*z]
\end{lstlisting}

\goodbreak
\noindent
Or coefficients in a fraction field:
\begin{lstlisting}[alsolanguage=C++]
/**/ use QQab ::= QQ[a,b];
/**/ K := NewFractionField(QQab);
/**/ use K[x,y,z];
/**/ I := ideal(x^3 -a*z, x^2*y^3 -b*y*z^2);
/**/ ReducedGBasis(I);
[x^3 -a*z,  x^2*y^3 -b*y*z^2,  y^3*z +(-b/a)*x*y*z^2]
\end{lstlisting}

One should note that in this last example $K$ is actually the field $\QQ(a,b)$
with no specialization of $a,b\in\QQ$.
So the \GBasis produced represents the
\textit{generic} case, meaning that every algebraic expression in $a,b$
which is not identically zero is considered to be non-zero.
The problem of considering all possible specializations of the
parameters is known as \textit{comprehensive \GBasis}, and is not
(yet) implemented in \cocoa.

Another family of computationally interesting rings in \cocoa is given by
\texttt{NewRingTwinFloat(BitPrec)}.
These will be presented in detail in Sections~\ref{sec:gin} and
\ref{sec:twinfloats}.

\section{Universal \GBases and Gr\"obner fans}
\label{sec:gfan}

There is a notion of \textbf{universal \GBasis} which is a \GBasis
for every term-ordering.
The \cocoa function \cocoacode{UniversalGBasis} will compute one
such basis; this function is based on the computation of
the Gr\"obner fan (a richer structure, described below) which gives all
possible reduced \GBases:
we can take the union of all of them to produce the universal basis.

The following example shows that the maximal minors of a $3 \times 4$
matrix of indeterminates form a
universal \GBasis of the ideal they generate:
\begin{lstlisting}[alsolanguage=C++]
/**/ use R ::= QQ[a,b,c,d,e,f,g,h,i,j,k,l];
/**/ I := ideal(minors(mat([[a,b,c,d],[e,f,g,h],[i,j,k,l]]),3));
/**/ indent(UniversalGBasis(I));
[
  d*g*j -c*h*j -d*f*k +b*h*k +c*f*l -b*g*l,
  d*g*i -c*h*i -d*e*k +a*h*k +c*e*l -a*g*l,
  d*f*i -b*h*i -d*e*j +a*h*j +b*e*l -a*f*l,
  c*f*i -b*g*i -c*e*j +a*g*j +b*e*k -a*f*k
]
/**/ EqSet(-1*gens(I), ReducedGBasis(I));
true
\end{lstlisting}




The \textbf{Gr\"obner fan} of an ideal was defined by Mora and Robbiano in 1988
(\cite{MorRob88}): it is a (finite) fan of polyhedral cones indexing
the reduced \GBases of the ideal.
This has been implemented by Jensen in his software \textit{Gfan} (\cite{Gfan})
which he has recently linked into \cocoa;
we note that \cocoa's fully general approach to representing
term-orderings was essential in making this integration possible.

The Gr\"obner fan is useful because several well-known
theoretical applications of \GBases rely on the existence of a
\GBasis of an ideal with prescribed properties, such as having
a certain cardinality, or comprising
polynomials of a specified degree, or all squarefree.
For example, if an ideal~$I\in K[x_1,\dots,x_n]$ has a \GBasis for
\textit{some} term-ordering comprising just quadrics, then the algebra
$K[x_1,\dots,x_n]/I$ is Koszul.

The function \texttt{GroebnerFanIdeals(I)} computes all reduced
\GBases of the ideal~$I$.  We have chosen to express the result as
a list of ideals, each generated by one of the various possible reduced \GBases:
the ideals are all ``the same'' but belong
to \textit{different} polynomial rings (remember that the term-ordering
is an intrinsic property of the polynomial ring).
An advantage of this approach is that
further computation with any of these ideals automatically takes
place in the corresponding polynomial ring equipped with
an appropriate term-ordering.
Furthermore, each of these ideals already knows its own reduced \GBasis,
whose value is thus immediately available (\ie~without any computation).


The following ideal, Example~3.9 from Sturmfels's book~\cite{Stu96},
has 360 distinct reduced \GBases:
\begin{lstlisting}[alsolanguage=C++]
/**/ use R ::= QQ[a,b,c];
/**/ I := ideal(a^5+b^3+c^2-1, a^2+b^2+c-1, a^6+b^5+c^3-1);
/**/ L := GroebnerFanIdeals(I);
/**/ len(L);
360
\end{lstlisting}

Since the computation easily becomes very cumbersome, it is
  interesting to see how it is progressing; for example, after setting the
  verbosity level to 10
(see Section~\ref{sec:verbosity}),
a \texttt{*} is printed every time a new
  Gr\"obner basis is added to the list (for more information see the
  manual by typing ``\texttt{?GroebnerFan}''):
\begin{lstlisting}[alsolanguage=C++]
/**/ use QQ[x,y,z];
/**/ I := ideal(x^3 +x*y -z, x^2 -y*z);
/**/ SetVerbosityLevel(10);
/**/ GF := GroebnerFanIdeals(I);
********
/**/ indent(GF);
\end{lstlisting}
\vskip -10pt
\begin{lstlisting}[alsolanguage=C++, basicstyle=\ttfamily\upshape\footnotesize]
[ ideal(x^2 -y*z,  x*y*z +x*y -z,  y^2*z^2 +y^2*z -x*z),
  ideal(x^2 -y*z,  x*z -y^2*z^2 -y^2*z,  y^3*z^2 +x*y +y^3*z -z),
  ideal(x^2 -y*z,  x*z -y^2*z^2 -y^2*z,  x*y +y^3*z^2 +y^3*z -z,
	  y^3*z^3 +2*y^3*z^2 +y^3*z -z^2),
  ideal(y*z -x^2,  x^3 +x*y -z),
  ideal(x*y +x^3 -z,  y*z -x^2,  x^3*z -z^2 +x^3),
  ideal(x*y -z +x^3,  y*z -x^2,  z^2 -x^3*z -x^3),
  ideal(z -x*y -x^3,  x*y^2 +x^3*y -x^2),
  ideal(z -x^3 -x*y,  x^3*y +x*y^2 -x^2) ]
\end{lstlisting}

Storing all the possible different (reduced)
\GBases is practicable only for small
examples; larger ideals may have thousands or even millions of different
\GBases.  Often we are interested only in those bases satisfying a
certain property.  So \cocoa
offers the function \texttt{CallOnGroebnerFanIdeals}
which calls a given function on each of the Gr\"obner fan ideals successively without
storing them all in a big list (which may not even fit in the computer's memory!).
Using this \cocoa function needs a little technical ability, but makes
it possible to tackle larger computations.

In the following example
we see explicitly that \cocoa represents some term-orderings
via matrices of integers.
Indeed, each such matrix is the only information necessary to
be able to recalculate the corresponding reduced Gr\"obner basis
(in this example, those having 3 elements).
See Section~\ref{sec:elim} for an example of how to ask \cocoa to
compute a \GBasis with a term-ordering given by a matrix.

\goodbreak
\begin{lstlisting}[alsolanguage=C++]
define PrintIfGBHasLen3(I)
 if len(GBasis(I))=3 then
   println OrdMat(RingOf(I));
   indent(ReducedGBasis(I));
 endif;
enddefine;

/**/ use R ::= QQ[a,b,c];
/**/ I := ideal(a^5+b^3+c^2-1, b^2+a^2+c-1, c^3+a^6+b^5-1);
/**/ CallOnGroebnerFanIdeals(I, PrintIfGBHasLen3);
\end{lstlisting}
\vskip -10pt
\begin{lstlisting}[alsolanguage=C++,basicstyle=\ttfamily\upshape\footnotesize]
matrix(ZZ,
 [[3, 7, 7],
  [3, 6, 8],
  [0, 0, -1]])
[b^2+c+a^2-1,
 a^5+c^2-b*c-a^2*b+b-1,
 c^3+b*c^2+2*a^2*b*c+a^4*b-a*c^2+a*b*c+a^3*b-2*b*c-2*a^2*b-a*b+b+a-1]
matrix(ZZ,
 [[6, 7, 14],
  [6, 5, 15],
  [0, 0, -1]])
[c+b^2+a^2-1,
 -b^6-3*a^2*b^4-3*a^4*b^2+b^5+3*b^4+6*a^2*b^2+3*a^4-3*b^2-3*a^2,
 a^5+b^4+2*a^2*b^2+a^4+b^3-2*b^2-2*a^2]
\end{lstlisting}

\section{Leading Term Ideals and ``gin''}\label{sec:gin}

 Let $P = K[x_1,\ldots, x_n]$ be a polynomial ring over a field
  $K$, and let $\sigma$ be a term-ordering on the power-products in
  $P$.  Let~$I$ be an ideal in $P$ then we define its \textbf{leading
    term ideal} with respect to $\sigma$, written $\LPP_\sigma(I)$, to
  be the ideal generated by the leading power-products of all non-zero
  polynomials
  in~$I$; some authors use the name ``initial ideal'' for this notion.
A generating set for $\LPP_\sigma(I)$ may easily be obtained:
we compute a reduced $\sigma$-\GBasis for~$I$ then collect the
$\sigma$-leading terms of the elements of the basis.
Remarkably $\LPP_\sigma(I)$ captures some interesting
  ``combinatorial'' information about the original polynomial ideal
$I$: for instance, its Hilbert series.  
Hence in CoCoA
calling \texttt{LT(I)} or \texttt{HilbertSeries(P/I)} actually
contains a ``hidden'' call to \texttt{GBasis(I)}.

\begin{lstlisting}[alsolanguage=C++]
/**/ use P ::= QQ[x,y,z];
/**/ I := ideal(y^20 -x^4*z^16,  x^12*z^3 -y^13*z^2);
/**/ LT(I);
ideal(y^13*z^2,  y^20,  x^12*y^7*z^3,  x^24*z^4)

/**/ HilbertSeries(P/I);
\end{lstlisting}
\vskip -10pt
\begin{lstlisting}[alsolanguage=C++,basicstyle=\ttfamily\upshape\footnotesize]
(1 + 2*t + 3*t^2 + 4*t^3 + 5*t^4 + 6*t^5 + 7*t^6 + 8*t^7 + 9*t^8 + 10*t^9
+ 11*t^10 + 12*t^11 + 13*t^12 + 14*t^13 + 15*t^14 + 15*t^15 + 15*t^16
+ 15*t^17 + 15*t^18 + 15*t^19 + 14*t^20 + 13*t^21 + 12*t^22 + 11*t^23
+ 10*t^24 + 9*t^25 + 8*t^26 + 7*t^27 + 6*t^28 + 5*t^29 + 4*t^30 + 3*t^31
+ 2*t^32 + t^33) / (1-t)
\end{lstlisting}


\medskip
A more sophisticated tool in Commutative Algebra is the
\textbf{generic initial ideal} of a polynomial ideal~$I$.
This is useful because it
encodes more geometrical properties of~$I$ into a monomial ideal.
It is defined as $\gin_\sigma(I) = \LPP_\sigma(\gamma(I))$
where $\gamma$ is a generic change of coordinates, \ie~$\gamma(x_j) =
\sum_{i=1}^n a_{ij}\,x_i$ in $K(a_{ij})[x_1,...,x_n]$.
And here we have to admit that the acronym \textit{gin} sounds nicer
than $gLT$!

The definition of \textit{gin} suggests an obvious algorithm for
computing it (see Section~\ref{sec:morecoefficients} for an example
with \textit{generic} coefficients).
However, 
even knowing that it is enough to consider a
\textit{triangular} change of coordinates
$\gamma(x_j) = \sum_{i=1}^j a_{ij}\,x_i$,
it quickly becomes apparent that the
coefficients in $K(a_{ij})$ grow to unwieldy sizes except for the very
simplest cases; so the obvious approach is utterly hopeless.  
Instead
we can pick an explicit, random change of coordinates $\tilde\gamma$, and then
compute $\LPP_\sigma(\tilde\gamma(I))$; the coordinate changes for which
$\gin_\sigma(I) = \LPP_\sigma(\tilde\gamma(I))$ form a Zariski-open set.
This approach can be used when~$K$ is infinite: if the
coefficients for the random change of coordinates are
chosen from a large set then $\LPP_\sigma(\tilde\gamma(I))$
will indeed be $\gin_\sigma(I)$  with high probability.
This is what \cocoa does.

While choosing random changes of coordinates with large coefficients
increases the probability of getting the correct result, it also
tends to produce large coefficients in the transformed polynomials.
In the example below the original polynomials have very small coefficients,
but there is a coefficient with almost 50 digits in the transformed
polynomials:
\begin{lstlisting}[alsolanguage=C++]
/**/ use P ::= QQ[x,y,z];
/**/ I := ideal(y^20 -x^5*z^6,  x^2*z^3 -y*z^2);
/**/ L := [sum([ random(-500,500)*indet(P,j) | j in 1..3])
                                                   | i in 1..3];
/**/ L;
[-414*x +341*y -141*z, -318*x +389*y +178*z, -498*x +498*y +28*z]
/**/ gamma := PolyAlgebraHom(P, P, L);
/**/ GI := ideal(apply(gamma, gens(I)));  GI;
\end{lstlisting}
\vskip -10pt
\begin{lstlisting}[alsolanguage=C++, basicstyle=\ttfamily\upshape\footnotesize]
ideal(111825899364055159629646472958188266490923110629376*x^20 +...,
      -21168433004832*x^5 +98376968843712*x^4*y - ...)
\end{lstlisting}
With coefficients like that,
computing the \GBasis of the transformed ideal over the rationals would be quite expensive!
Thus, one needs to strike a balance between picking coefficients from a wide
range, so the transformation is ``generic enough'', but not
so wide that there is excessive growth in the coefficients of transformed ideal generators.

\medskip
To avoid the costs of computing with large coefficients,
the implementation for computing $\gin$ in \cocoa uses a special
representation for rational coefficients,
namely \textit{twin-floats} (see Section~\ref{sec:twinfloats}).
The \GBasis of the twin-float transformed ideal will have only approximate
twin-float coefficients, but this does not matter because we need only
the leading power-products of the polynomials in the basis.

Twin-float numbers have fixed-precision (so do not grow in size the
way rational numbers do), and employ heuristics to verify the
correctness of results.  This allows the implementation to make random
coefficient choices from a wide range (in fact, integers between
$-10^6$ and $10^6$) without paying the price for calculating with
transformed polynomials having cumbersome rational coefficients.
If the initially chosen precision for the twin-floats is too low, this will be
signalled; and the computation will be automatically restarted with a higher
precision.

\cocoa's function \texttt{gin} does all this behind the scenes.
Moreover, it also tries a second random change of coordinates,
just to make sure it gets the same leading term ideal.
If the results differ,
\cocoa keeps trying further random changes of coordinates until it gets
the same answer twice in succession~---~though we have never seen this
happen when picking random coefficients in the range $(-10^6,10^6)$.
The internal workings can be seen via printed messages
with the appropriate verbosity level (see Section~\ref{sec:verbosity}).
\goodbreak
\begin{lstlisting}[alsolanguage=C++]
/**/ SetVerbosityLevel(50);
/**/ J := gin(I);
RandIdeal: change coord = [
  -7426*x,
  695955*x +168758*y,
  -239080*x +304634*y +480790*z
]
TryPrecisions: -- trying with FloatPrecision 64
TryPrecisions: -- trying with FloatPrecision 128
RandIdeal: change coord = [
  -499447*x,
  -732749*x -840921*y,
  -466314*x -691911*y +554086*z
]
TryPrecisions: -- trying with FloatPrecision 128

/**/ J;
ideal(x^5,  x^4*y^16,  x^3*y^18,  x^2*y^20,  x*y^22,  y^24)
\end{lstlisting}

Since the $\gin$ ideal with respect to
the 
 ordering
\texttt{StdDegRevLex}
 has many interesting properties,
there is also the function \texttt{rgin}
which computes it, independently of the term ordering inherent in the
polynomial ring.

\section{Elimination and related functions}
\label{sec:elim}

Elimination means:
given an ideal $I \in K[t_1,\dots,t_s, \; x_1,\dots,x_n]$,
find a set of generators of the ideal $I \cap K[x_1,\dots,x_n]$
where the indeterminates $\{t_1,\dots,t_s\}$ have been ``eliminated''.
Elimination is a central topic in Computational Commutative Algebra
(see for example the text book
by Kreuzer and Robbiano~\cite{KR00}, Sec.~3.4) and its
applications are countless.

Given its usefulness, elimination is an operation offered in
  almost
all Computer Algebra Systems.   
In general, such elimination
functions internally compute a \GBasis with respect to an
\textbf{elimination ordering} for the subset of indeterminates to be eliminated: 
with such an ordering the subset of polynomials in the \GBasis whose leading terms
are not divisible by any of the $t_j$ are exactly the generators we seek for
the ideal $I\cap K[x_1,\dots,x_n]$.

In the example below we see the process we described and compare it
with the actual output of \cocoa's own function \texttt{elim}.
Note that in both cases the generators are not minimal, but they are
indeed a \GBasis of the elimination ideal (wrt~the restriction of
the elimination term-ordering used).

\begin{lstlisting}[alsolanguage=C++]
/**/ M := ElimMat([1], 4);  M;
matrix(ZZ,
 [[1, 0, 0, 0],
  [1, 1, 1, 1],
  [0, 0, 0, -1],
  [0, 0, -1, 0]])

/**/ P := NewPolyRing(QQ, "t, x,y,z", M, 0);  // 0: no grading
/**/ use P;
/**/ I := ideal(x-t, y-t^2, z-t^3);
/**/ GBasis(I);
[t -x,  x^2 -y,  x*y -z,  y^2 -x*z]
/**/ elim([t], I);
ideal(x^2 -y,  x*y -z,  y^2 -x*z)
/**/ MinSubsetOfGens(ideal(x^2 -y,  x*y -z,  y^2 -x*z));
[x^2 -y,  x*y -z]
\end{lstlisting}

The simple example above shows a particular application of
\texttt{elim}: finding the \textbf{presentation of an algebra}
$K[f_1,...,f_n]\simeq K[x_1, ..., x_n]{/}I$.  More precisely,
let   $f_1, ..., f_n \in K[t_1, \dots, t_s]$,
where $\{t_1, \dots, t_s\}$ is another set of indeterminates
(viewed as parameters) and consider the $K$-algebra homomorphism
$$
\phi: K[x_1, \dots, x_n] \longrightarrow K[t_1, \dots, t_s]
 \text{ given by }
x_i \mapsto f_i \hbox{\rm \quad  for\ }  i=1, \dots, n
$$
Its kernel 
 is a prime ideal;
 the general problem of \textbf{implicitization} 
(for a polynomial parametrization) is to
find a set of generators for this ideal.

The \GBasis elimination technique consists of defining the ideal $J =
\ideal{ x_1-f_1, \dots, x_n-f_n }$ in the ring $K[ t_1, \dots, t_s,
x_1, \dots, x_n]$ and \textit{eliminating} all the parameters $t_i$,
as we saw in the example.

Unfortunately this extraordinarily elegant tool often turns out to be quite
inefficient, resulting in long and costly computations.
Knowing how to exploit special properties of
a given class of examples
might make a huge difference.

\subsection{Toric}

If the algebra we want to present is generated by power-products then
the elimination can be computed by the \cocoa function
 \cocoacode{toric}; toric ideals are prime and generated by binomials.

We consider $\QQ[t,t^2,t^3]$ again, and also $\QQ[t^3,t^4,t^5]$:
\begin{lstlisting}[alsolanguage=C++]
/**/ use QQ[x,y,z];
/**/ toric(RowMat([1,2,3])); // just the list of exponents
ideal(-x^2 +y, x^3 -z)

/**/ use QQ[x,y,z];
/**/ toric(RowMat([3,4,5]));
ideal(y^2 -x*z,  x^3 -y*z,  x^2*y -z^2)
\end{lstlisting}

With a very slightly more challenging example we can clearly measure
the advantage in using the specialized function \cocoacode{toric}
over the general function \cocoacode{elim}:
\begin{lstlisting}[alsolanguage=C++]
/**/ use R ::= ZZ/(2)[x[1..6], s,t,u,v];
/**/ L := [s*u^20, s*u^30, s*t^20*v, t*v^20, s*t*u*v, s*t^2*u];
/**/ ExpL := mat([[ 1,   1,   1,   0,   1,   1],
                  [ 0,   0,  20,   1,   1,   2],
                  [20,  30,   0,   0,   1,   1],
                  [ 0,   0,   1,  20,   1,   0]]);

/**/ I := ideal([x[i] - L[i] | i in 1..6]);
/**/ t0 := CpuTime();  IE := elim([s,t,u,v], I);  TimeFrom(t0);
9.274
/**/ t0 := CpuTime();  IT := toric(ExpL);  TimeFrom(t0);
0.032
\end{lstlisting}

The \cocoa function \cocoacode{toric} employs a non-deterministic
algorithm: so the actual set of ideal generators produced might vary.

For further details on the algorithms implemented in \cocoa see
Bigatti, La~Scala, Robbiano \cite{BigLScRob99}.
That article
describes three different algorithms; the default one in \cocoa is
EATI (\textit{Elimination Algorithm for Toric Ideals}).

For more details on the specific function \cocoacode{toric} type
\texttt{?toric} into \cocoa (or read the PDF manual, or the html manual
on the web-site).

\subsection{Implicitization of hypersurfaces}
\label{sec:implicit}

As mentioned earlier elimination provides a general solution to the
implicitization problem, but this solution is more elegant than
practical.  We can do rather better in the special case of
implicitization of a hypersurface.
One immediate feature is that the
result is just a single polynomial since the eliminated ideal
must be principal.

It is well-known that Buchberger's algorithm usually works
better with homogeneous ideals (though there are sporadic exceptions).
Yet the very construction of the eliminating ideal
$J = \ideal{ x_1-f_1, \dots, x_n-f_n }$
looks intrinsically non-homogeneous.

But with a little, well-guided effort we can transform the problem into
the calculation of a \GBasis of a homogeneous ideal:
we take a new indeterminate (say~$h$)
and use it to homogenize each $f_j$ to produce $F_j$.
Now we can work with the ideal 
$J' = \ideal{ x_1-F_1, \dots, x_n-F_n }$
made homogeneous by giving weights to the $x_i$ indeterminates: we set
$\deg(x_i) = \deg(f_i)$ for each~$i$.

Since we are in the special case of a hypersurface, it can be shown
that the (non-zero) polynomial of lowest degree in 
$J' \cap K[x_1,\ldots,x_n,h]$ is unique up to scalar
multiples; its dehomogenization is then the polynomial we seek!
We get two advantages from the homogeneous ideal $J'$: we gain efficiency
by using Buchberger's algorithm degree-by-degree, and we can stop as
soon as the first basis polynomial is found~---~most probably there
will still be many pairs to process.
See Abbott, Bigatti, Robbiano~\cite{AbbBigRob17} for all details and
proofs, and also how to ``correctly homogenize'' parametrizations 
defined by rational functions.

\begin{lstlisting}[alsolanguage=C++]
/**/ use P ::= QQ[s,t, x,y,z];
/**/ elim([s,t],  ideal(x-s^2, y-s*t, z-t^2) );
ideal(y^2 -x*z)

/**/ use R ::= QQ[s,t];
/**/ P ::= QQ[x,y,z];
/**/ ImplicitHypersurface(P, [s^2, s*t, t^2], "ElimTH");
ideal(y^2 -x*z)
\end{lstlisting}

In the same paper we describe another algorithm which uses a
completely different technique, a variant of
the Buchberger-M\"oller algorithm (see Section~\ref{sec:BM}),
based on linear algebra.  It is well-suited to low
degree 
hypersurfaces.
\begin{lstlisting}[alsolanguage=C++]
/**/ ImplicitHypersurface(P, [s^2, s*t, t^2], "Direct");
ideal(y^2 -x*z)
\end{lstlisting}

For the case of rational
coefficients, we use a modular approach in both algorithms:
we compute the result modulo several
primes, combine these using Chinese Remaindering, and finally reconstruct the
rational coefficients of the answer using the fault-tolerant rational
reconstruction described in Section~\ref{sec:faulttolerant}.

\subsection{MinPoly}\label{sec:minpoly}

Another popular application of elimination is
for finding univariate polynomials in an ideal.
If $I\subset P = K[x_1,\dots,x_n]$ is a  0-dimensional ideal, then
we know that $I\cap K[x_i]$ is a principal ideal generated by some
univariate polynomial
$g_i(x_i)\ne0$ which can be obtained by eliminating all $x_j$ with $j \neq i$.
These polynomials are used in several operations, such as
computing the radical of a zero-dimensional ideal, or solving
polynomial systems (see~\cite{KR3} and \cite{AbbBigPalRob17}).

This idea generalizes in a natural way to the following problem.  Let $I\subset P$ be a
0-dimensional ideal, and let $f$ be any polynomial in $P$, find
$\minpoly{f,I}(z)\in K[z]$, the
minimal polynomial of $\bar f\in P/I$, or equivalently,
the univariate polynomial of minimum
degree whose evaluation at $f$ yields an element of $I$.
The corresponding algorithms have been recently implemented in CoCoA:

\begin{lstlisting}[alsolanguage=C++]
/**/ use P ::= QQ[x,y,z];
/**/ L := [ x^2-z^2, (y-3)*(y+2)*(y^3-2), z^3-1];
/**/ I := ideal(L);
/**/ IsZeroDim(I);
true

/**/ MinPolyQuot(x, I, x); -- 3rd arg is the indet for the answer
x^6 -1

/**/ f := x -2*y +3*z;
/**/ t0 := CpuTime();  MP := MinPolyQuot(f, I, x);  TimeFrom(t0);
0.036
/**/ MP;
\end{lstlisting}
\vskip -10pt
\begin{lstlisting}[alsolanguage=C++,basicstyle=\ttfamily\upshape\tiny]
x^30 +12*x^29 -84*x^28 -1544*x^27 +384*x^26 +62688*x^25 +119168*x^24 -629760*x^23 -4664832*x^22
-33803264*x^21 +107753472*x^20 +1318662144*x^19 -3480064000*x^18 -20059865088*x^17
+151993466880*x^16 -50058002432*x^15 -1931977162752*x^14 +9312278544384*x^13 +1002303913984*x^12 
-113944836440064*x^11 +553708192530432*x^10 +720752546414592*x^9 -6749908862238720*x^8
+4995175176732672*x^7 +33972228030726144*x^6 +22154393721765888*x^5 -21399162914340864*x^4
-112685231584051200*x^3 +3245139849904128*x^2 -3103199770705920*x -16498446852685824
\end{lstlisting}
Needless to say, even in this small example, the standard elimination
approach is considerably slower:
\begin{lstlisting}[alsolanguage=C++]
/**/ use Paux ::= QQ[x,y,z, aux];
/**/ phi := PolyAlgebraHom(P, Paux, [x,y,z]);
/**/ J := ideal(apply(phi,L)) + ideal(aux - phi(f));
/**/ t0 := CpuTime();  JE := elim([x,y,z], J);  TimeFrom(t0);
1.850
\end{lstlisting}

As we did for hypersurface implicitization (in Section~\ref{sec:implicit}),
when computing with rational coefficients
we use a modular approach
and the fault-tolerant rational
reconstruction described in Section~\ref{sec:faulttolerant}.

The good computational speed of \texttt{MinPolyQuot} is the key point for a new
algorithm for computing the primary decomposition of zero-dimensional ideals.
See~\cite{AbbBigPalRob17} for details on
the algorithms for \texttt{MinPolyQuot} and some interesting applications.

\goodbreak
\begin{lstlisting}[alsolanguage=C++]
/**/ PD := PrimaryDecomposition0(I);
/**/ indent([IdealOfGBasis(Qi) | Qi in PD]);
\end{lstlisting}
\vskip -10pt
\begin{lstlisting}[alsolanguage=C++,basicstyle=\ttfamily\upshape\tiny]
[
  ideal(y +2,  x +1,  z -1),
  ideal(y -3,  x +1,  z -1),
  ideal(x +1,  z -1,  y^3 -2),
  ideal(y +2,  x -1,  z -1),
  ideal(y -3,  x -1,  z -1),
  ideal(x -1,  z -1,  y^3 -2),
  ideal(y +2,  z^2 +z +1,  x +z),
  ideal(y -3,  z^2 +z +1,  x +z),
  ideal(z^2 -x +1,  y^3 -2,  x +z),
  ideal(y +2,  z^2 +z +1,  x -z),
  ideal(y -3,  z^2 +z +1,  x -z),
  ideal(z^2 +x +1,  y^3 -2,  x -z)
]
\end{lstlisting}


\section{Ideals of Points, 0-Dimensional Schemes}\label{sec:BM}

Let $X$ be a non-empty, finite set of points in $K^n$, then
the set of all polynomials in $K[x_1,\ldots,x_n]$ which vanish
at all points in~$X$ is an ideal, $I_X$.  One reason this ideal
is interesting is because it captures the
``ambiguity'' present in a polynomial function which has been
interpolated from its values at the points of $X$.  How best to
compute a set of generators for $I_X$, or a \GBasis, knowing just the
points~$X$?

If $X$ contains a single point $(a_1, \ldots, a_n)$ then we can
write down immediately a \GBasis, namely $[x_1-a_1, \ldots, x_n-a_n]$.
If $X$ contains several points we could just intersect the ideals
for each single point, and these intersections may be
determined via \GBasis computations; while fully effective and mathematically
elegant this approach is computationally disappointing.

A far more efficient method is the Buchberger-M\"oller algorithm~\cite{BM82}.
Somewhat astonishingly it uses just simple linear algebra to determine
the \GBasis.  In~\cite{AbbBigKreRob00} there is a detailed complexity analysis of
the original algorithm,
and also an extension to the projective case.
It was later further generalized to zero-dimensional
schemes~\cite{AbbKreRob05}, where it turned out that it also incorporates
the well-known FGLM algorithm for ``changing term-ordering'' of a
\GBasis (see \cite{FGLM}).

Much as we have seen in the previous sections for computing with rational coefficients,
the Buchberger-M\"oller algorithm also benefits from a modular approach,
and naturally the \cocoa implementation uses this technique.

\begin{lstlisting}[alsolanguage=C++]
/**/ P ::= QQ[x,y];
/**/ points := mat([[10, 0], [-10, 0], [0, 10], [0, -10],
                     [7, 7], [-7, -7], [7, -7], [-7, 7]]);
/**/ indent(IdealOfPoints(P, points));
ideal(
  x^2*y +(49/51)*y^3 +(-4900/51)*y,
  x^3 +(51/49)*x*y^2 -100*x,
  y^4 +(-2499/2)*x^2 +(-2699/2)*y^2 +124950,
  x*y^3 -49*x*y  )
\end{lstlisting}

The use of simple linear algebra in the Buchberger-M\"oller
 algorithm makes it a good
candidate for identifying ``almost-vanishing'' polynomials for sets of
\textit{approximate points:} for instance, the points in the example
above ``almost
lie on'' a circle of radius $9.95$ centred on the origin, though we cannot tell this from the \textit{exact} \GBasis.

In fact, the notion of
\GBasis does not generalize well to an ``approximate context'' because the
algebraic structure of a \GBasis is determined by Zariski-closed conditions
(\ie~the structure is valid when certain polynomials vanish); instead, the
notion of a \textit{Border Basis} is better suited since the validity of its structure
depends on a Zariski-open condition (\ie~provided a certain
polynomial \textit{does not} vanish).  So long as the approximate points are not
too few nor too imprecise the NBM (\textit{Numerical Buchberger-M\"oller})
algorithm can compute at least a partial
Border Basis, and this should identify any ``approximate polynomial
conditions'' which the points the almost satisfy
(see Abbott, Fassino, Torrente \cite{AbbFasTor08} and 
Fassino \cite{Fas10}).
We can ask \cocoa to
allow a certain approximation on the coordinates of
the points:
\begin{lstlisting}[alsolanguage=C++]
/**/ epsilon := [0.1, 0.1]; // coord approximation 0.1
/**/ AP01 := ApproxPointsNBM(P, mat(points), mat([epsilon]));
/**/ indent(AP01.AlmostVanishing);
[
  x^2 +(4999/5001)*y^2 -165000/1667, // almost a circle
  x*y^3 -49*x*y,
  y^5 -149*y^3 +4900*y
]

/**/ epsilon := [0.01, 0.01]; // coord approximation 0.01
/**/ AP001 := ApproxPointsNBM(P, mat(points), mat([epsilon]));
/**/ indent(AP001.AlmostVanishing);  // not "epsilon-near" a conic
[
  x^2*y +(49/51)*y^3 +(-4900/51)*y,
  x^3 +(51/49)*x*y^2 -100*x,
  y^4 +(-2499/2)*x^2 +(-2699/2)*y^2 +124950,
  x*y^3 -49*x*y
]
\end{lstlisting}

\section{\GBases and rational coefficients}\label{sec:rat}

It is well known that computations with coefficients in $\QQ$ can
often be very costly in terms of both time and space.  For \GBases
over~$\QQ$ we are free to multiply the polynomials by any non-zero
rational; so we can clear denominators and remove integer content.
Avoiding rational arithmetic this way
 does yield some benefit, but is not wholly satisfactory.

Sometimes the \GBasis has complicated coefficients (\ie~we mean \textit{big
numerators and denominators}), but more often the coefficients in the
answer are reasonably sized, while the computation to obtain them
involved far more complicated coefficients: this problem is known
as \textit{intermediate coefficient swell}.

The phenomenon of coefficient swell is endemic in computer algebra,
and many techniques have been investigated to tackle this problem.  We
illustrate 
two techniques used in \cocoa.

\subsection{TwinFloat}\label{sec:twinfloats}

\cocoa offers floating-point arithmetic with a heuristic verification of
correctness: the aim is to combine the speed of floating-point
computation with the reliability of exact rational arithmetic~---~for
a fuller description see the article~\cite{Abb12}.  Normally a twin-float
computation will produce either a good approximation to the correct
result or an indication of failure; strictly, there is a very small
chance of getting a wrong result, but this never happens in practice.

To perform a computation with twin-floats the user must first specify
the required precision; \cocoa will then perform the computation checking
heuristically that the
result of every twin-float operation has at least that precision.  If the check
fails then \cocoa signals an ``insufficient precision'' error; the
user may then restart the computation specifying a higher precision.
Although twin-float values are, by definition, approximate, all input
values are assumed to be exact (so they can be converted to a twin-float
of any precision).

It is also possible to convert a twin-float value to an exact
rational number.  Like all other twin-float operations, this
conversion may fail because of ``insufficient precision''.  Printing
out a twin-float value automatically attempts conversion to
a rational as rationals are easier to read and comprehend
in the context of exact computations.

  \begin{lstlisting}[alsolanguage=C++]
/**/ RR16 := NewRingTwinFloat(16);
/**/ use RR16_X ::= RR16[x,y,z];
/**/ f := 12345678*x+1/456789;
/**/ f;          // both coeffs are printed as rationals
12345678*x +1/456789
/**/ f * 10^3;   // first coeff is printed in "floating-point"
0.12345678*10^11*x +1000/456789
/**/ f * 10^5;   // both coeffs are printed in "floating-point"
0.12345678*10^13*x +0.2189194573
\end{lstlisting}

  Twin-floats include a (heuristically verified) test for zero; this
  means it is possible to compute \GBases with twin-float coefficients.
  One reason for wanting to do this is that often the computation of a \GBasis over the
  rationals involves ``complicated fractions'' (\ie~whose numerator and
  denominator have many digits), and arithmetic with such complicated
  values can quickly become very costly.  In contrast, with twin-floats
  the arithmetic has fixed cost (dependent on the precision chosen, of
  course).
  These characteristics are exploited in \cocoa for the computation of
  the $\gin$,  described in Section~\ref{sec:gin}.

\subsection{(Fault-tolerant) Rational reconstruction}
\label{sec:faulttolerant}

A widely used technique for avoiding intermediate coefficient swell is
to perform the computation modulo one or more prime numbers, and then
lift/reconstruct the final result over $\QQ$.  We call this the
\textbf{modular approach}.  There are two general classes of method:
\textbf{Hensel Lifting} and \textbf{Chinese Remaindering}, the first
is not universally applicable but does work well for polynomial gcd
and factorization, while the second is widely applicable and works
well in most other contexts.

The modular approach has been successfully used in numerous contexts,
here are a few examples: polynomial factorization~\cite{Zas69},
determinant of integer matrices~\cite{AbbBroMul99}, ideals of points
(see Section~\ref{sec:BM}), implicitization (see
Section~\ref{sec:implicit}), 
and minimal polynomial (see Section~\ref{sec:minpoly}).

In any specific application
there are two important aspects which must be addressed before a
modular approach can be adopted, and there is no universal technique
for addressing these issues:
\begin{itemize}[parsep=-0pt]
\item
  knowing how many different primes to consider to guarantee the result
  (\ie~find a realistic bound for the size of coefficients in the answer);
\item
handling \textit{bad primes}: namely those whose related computation follows a
different route, yielding an answer with the wrong ``shape''
(\ie~which is not the modular reduction of the correct, non-modular result).
\end{itemize}

In the context of \GBases we do not have good, general solutions to either
of these issues.  One of the first successes in applying modular
techniques to \GBasis computation appeared in~\cite{Arn03}.  Finding
good ways to employ a modular approach for \GBases is still an active area.
\cocoa does not currently use a modular approach for general \GBasis
computations.

A vital complement to the modular computation is the reconstruction of
the final, rational answer from the modular images.  \cocoa offers
functions for
\begin{itemize}
\item combining two residue-modulus pairs into a single
``combined'' residue-modulus pair (using the Chinese Remainder
Theorem);
\item determining a ``simple''
rational number corresponding to a residue-modulus pair; this is called
\textbf{rational reconstruction}.  
\end{itemize}
Correct reconstruction can still be
achieved even in the presence of a few ``faulty residues''
(see~\cite{Abb16}); this fault-tolerance was exploited in the
functions for hypersurface implicitization (see Section~\ref{sec:implicit}).

Here we see how two modular images can be combined in \cocoa (using
\cocoacode{CRTPoly}), and then the correct rational
result is reconstructed from the combined residue-modulus pair
(using \cocoacode{RatReconstructPoly}).
\begin{lstlisting}[alsolanguage=C++]
/**/ P1 ::= ZZ/(12347)[x];
/**/ P2 ::= ZZ/(23459)[x];
/**/ ReadExpr(P1, "x/1234-1/5");  // modular image in P1
-5293*x -4939
/**/ ReadExpr(P2, "x/1234-1/5");  // modular image in P2
-1806*x -4692

/**/ use P ::= QQ[x];
/**/ combined := CRTPoly(-5293*x-4939, 12347,
                         -1806*x-4692, 23459);
/**/ combined;  // in P
record[modulus := 289648273, residue := 79571122*x +115859309]
/**/ RatReconstructPoly(combined.residue, combined.modulus);
(1/1234)*x -1/5
\end{lstlisting}

\section{Gr\"obner bases in C++ with CoCoALib}\label{sec:cocoalib}

As mentioned in Section~\ref{sec:cocoalib-vs-cocoa5}, our aim is
to make computation using CoCoALib as easy as using \cocoav.  To
illustrate this,
here is the first example from Section~\ref{sec:GBwithEase} but in C++:
\begin{lstlisting}[alsolanguage=C++]
ring P = NewPolyRing(RingQQ(), symbols("x,y,z"));
ideal I = ideal(ReadExpr(P, "x^3 + x*y^2 - 2*z"),
                ReadExpr(P, "x^2*y^3 -y*z^2")  );
cout << GBasis(I);
\end{lstlisting}
In comparison to \cocoav, this C++ code is more cumbersome and involved,
though we maintain that it is still reasonably comprehensible (once
you know that \verb|cout <<| is the C++ command for printing).


We have designed \cocoav and \cocoalib together with the aim of
making it easy to develop a prototype implementation in \cocoav,
and then convert the code into C++.
To facilitate this conversion we have, whenever
possible, used the same function names in both \cocoav and \cocoalib,
and we have preferred traditional ``functional'' syntax in CoCoALib over
object oriented ``method dispatch'' syntax
(\textit{e.g.}~\texttt{GBasis(I)} rather than~\texttt{I.GBasis()}).
This means that most of the \cocoav examples given here require only
minor changes to become equivalent C++ code for use with \cocoalib.

To maintain the ``friendly'' tradition of \cocoa software for
mathematicians, and to extend it to ``mathematical programmers'', our design
of \cocoalib follows these aims:
\begin{itemize}[itemsep=0pt]
\item Designed to be \textbf{easy and natural to use}
\item  Motto: \textbf{``No nasty surprises''} (avoid ambiguities)
\item Execution speed is \textbf{good}
\item Well-documented, including \textbf{many example programs}
  \smallskip
\item \textbf{Free and open source} C++ code (GPL3 licence)
\item Source code is \textbf{clean and portable} (currently C++03)
\item Design respects the underlying \textbf{mathematical} structures\\
(using C++ inheritance, no templates)
\item \textbf{Robust} exception-safe, thread-safe
\end{itemize}

\section{Conclusion}

The \cocoa software aims to make it easy for everyone to
use \GBases, whether directly or indirectly through some other
function.  The \cocoav system is designed to be welcoming to
those with little computer programming experience, while
the \cocoalib library aims to make it easy for experienced
programmers to use \GBases in their own programs.

We hope this helps everyone to have their \GBasis!

\goodbreak

\end{document}